\documentclass[12pt,amstex]{article}
\usepackage[usenames]{color}

\usepackage{amscd}
\usepackage{amsmath}
\usepackage{amssymb}
\usepackage{graphicx}

\newtheorem {theo} {\bf Theorem} [section]
\newtheorem {prop} [theo] {\bf Proposition}

\newtheorem {lem} [theo] {\bf Lemma}
\newtheorem {defn} [theo] {\bf Definition}

\newtheorem {rem} [theo] {\bf Remark}
\newcommand{\QED}{\hfill \CaixaPreta \vspace{6mm}}
\def\CaixaPreta{\vrule Depth0pt height6pt width6pt}
\newenvironment{remark}{\begin{rem}\begin{rm}}{\end{rm}
\nopagebreak\hfill$\Box$\end{rem}\par\vspace{\partopsep}}

\newcommand{\qed}{\nopagebreak\hfill{\vrule width6pt height6pt depth0pt}}

\newcommand{\be}{\begin{eqnarray}}
\newcommand{\ee}{\end{eqnarray}}
\newcommand{\benn}{\begin{eqnarray*}}
\newcommand{\eenn}{\end{eqnarray*}}
\newcommand{\bse}{\begin{equation}}
\newcommand{\ese}{\end{equation}}
\newcommand{\bsenn}{\begin{displaymath}}
\newcommand{\esenn}{\end{displaymath}}
\newcommand{\logand}{\;\;{\rm and }\;\;}
\newcommand{\for}{\;\;{\rm for }\;\;}

\newcommand{\where}{\;\;{\rm where }\;\;}
\newcommand{\with}{\;\;{\rm with }\;\;}

\newcommand{\R}{\mathbb{R}}

\newtheorem{corollary}{Corollary}

\newtheorem{cond}[theo]{\bf Condition}
\newcommand{\boldm}[1]{\mbox{\boldm\alphath$#1$}}

\newcommand{\half}{\frac{1}{2}}
\newcommand{\smhalf}{ {\scriptstyle{\half}} }

\begin{document}

\title{The Bunching and Monotonicity Properties of Families of Probability Distributions}
\author{ \\ 
S. Portnoy\thanks{Stat. Depts., Portland State University and University of Illinois; e-mail: spor2@psu.edu},\,  
N. Torrado\thanks{Math. Dept., Universidad Aut\'{o}noma de Madrid and Instituto de Ciencias Matm\'{a}ticas, 28049 Madrid; 
e-mail: nuria.torrado@uam.es; nuria.torrado@icmat.es} 
\, and J.J.P. Veerman\thanks{Maseeh Dept. of Math. and Stat., Portland State University,
	Portland, OR, USA; e-mail: veerman@pdx.edu.}}
\maketitle

\bigskip

\begin{abstract}
	
Measuring the concentration of random variables is a fundamental concept in probability and statistics.
Here, we explore a type of concentration measure for continuous random variables with bounded support 
and use it to provide a notion of stochastic order by concentration. We give an 
application  to the Beta family of distributions, and specifically to the
one-parameter subfamily with constant mean. This leads to using U.S. household income data to
fit generalized Beta distributions and offer a new measure of income concentration.

\end{abstract}

\bigskip

{\bf Keywords:} Stochastic Order, Concentration, Monotonicity, Beta Distribution

\medskip

{\bf MSC classification:}  60E15, 62E99

\bigskip



\section{Motivation and Introduction}
 \label{chap:intro}
\setcounter{figure}{0} \setcounter{equation}{0}

A fundamental focus in probability theory lies in examining the degree of concentration of a random variable. 
This concept involves investigating how closely or loosely the values of a random variable cluster
around a specific value, often its central measure, such as the mean or median.
The degree of concentration provides insights into the variability and predictability of the variable's outcomes. 
Understanding concentration is essential in various fields, including statistics, economics, and risk analysis,
as it enables us to assess the reliability and stability of the random variable in question.

Here we refer to the measure of concentration as ``bunching".  It was introduced to solve a motivating problem
concerning a one-parameter sub-family of Beta distributions (see the Conclusions Section for details). The measure
was developed some time ago (see \cite{Shak1}  and references there). Our basic theoretical result (Lemma 3.3)
provides a sufficient condition for bunching that is essentially equivalent to the result of Corollary 2.1 in 
\cite{Ramos}. The new contributions here are the following: a self-contained proof of the basic result, application to a one-parameter sub-family of the Beta family, and a specific application to income distribution suggesting the use of the Beta parameter as an alternative measure of inequality (with comparison to the Gini index).

There are several ways to measure the concentration of random variables. One of them is through concentration inequalities, 
which provide valuable bounds on how much random variables deviate from a specific value, typically the expected value. 
These inequalities are instrumental in helping us to measure the degree of concentration within a given distribution and 
understand how well the data points cluster around a specific value. They play a crucial role in diverse fields, from 
probability theory to statistics and machine learning, by allowing us to make informed assessments of the variability
and stability of random variables and their associated distributions. In particular, concentration inequalities can be
divided into two major groups: those that are distribution-free and those that are dependent on a specific distribution.
In the first group, we encounter well-known inequalities such as Chebyshev's inequality, Markov's inequality,
Chernoff's inequality, Hoeffding's inequality, and Bernstein's inequality, to name a few. Chebyshev's inequality
provides an upper bound for the probability that a random variable deviates by more than a certain number of
standard deviations from its mean. Markov's inequality, on the other hand, offers an exonential upper bound and
focuses on the probability that a non-negative random variable is greater than or equal to a specific value.
In contrast, Chernoff's and Hoeffding's inequalities are designed for independent random variables, aiming to
quantify the probability that the sum of these random variables deviates from its mean.
Bernstein's inequality serves as a generalization of Hoeffding's inequality and provides bounds for the probability
that the sum of independent random variables significantly deviates from its mean.

In the second group, Chvátal's conjecture focussed on the Binomial distribution and 
has recently attracted the attention of several researchers.
Specifically, Chvátal conjectured that for any given $n$, the probability of a binomial random variable $B(n, m/n)$
with integer expectation $m$ is smallest when $m$ is the integer closest to $2n/3$. \cite{Jason}  showed that this
holds when $n$ is large. \cite{Sun} proved that the sequence $q_k = 1 - p_k = P\{B(n, kn) \geq k+1\}$ is
strictly unimodal with the mode $k_0$ being the integer closest to $2n/3$ for any fixed $n \geq 2$. \cite{Barabesi}
established that Chvátal's conjecture is indeed true for every $n \geq 2$. Motivated by these works, \cite{Li} studied
the cases of Poisson, geometric, and Pascal distributions.

Another way to measure the concentration of a random variable is to examine whether the probability that the random variable
is less than or equal to a specific value exhibits a monotonic behavior. Apparently, the first work discussing such monotonicity
properties was \cite{Ghosh}, which provided results for the Chi-square, Fisher-Snedecor F, and Student's t-distributions.
Subsequently, \cite{Alam} explored the case of the incomplete Beta function. A recent work studying the case of the
gamma distribution is \cite{Pinelis}.
Furthermore, the concentration of two random variables can be compared using certain stochastic orders, 
such as the convex, dispersive, and right spread orders. To explore these and other stochastic orders,
refer to Shaked and Shanthikumar's monograph \cite{Shak}. 
Some works where characterizations and applications of these orders are studied include the following.
\cite{Belzunce} characterized the right spread order through the increasing convex order.
\cite{Ma} investigated convex orders for linear combinations of random variables.
In \cite{Komisarski}, they study the convex order between convolution polynomials of finite Borel measures.
\cite{casta} explored the increasing concave orderings of linear combinations of ordered statistics, applying them
to issues related to social well-being.
In \cite{che}, they examine two categories of optimal insurance decision problems associated with the convex order,
where the objective function or the premium valuation is a general function of the expected value, Value-at-Risk,
and Average Value-at-Risk of the loss variables.

In this work, we investigate a type of concentration measure for continuous random variables with bounded support
within the interval $[0,1]$ which we refer to as the bunching property. Specifically, we demonstrate that under certain
conditions, there exists a unique point $x^{*}$ around which the distribution is more bunched (or concentrated).
Additionally, we study a continuity and monotonicity property of such a point.
The structure of this article is as follows. Section 2 presents the definitions that will be employed throughout the manuscript.
Section 3 is devoted to the main results for continuous random variables with bounded support within the interval $[0,1]$.
Theorem \ref{thm:bunching} provides sufficient conditions
under which one member of a  one-parameter family of distributions
is more bunched than another. Section 4 considers the Beta family,
$B(\alpha, \, \beta)$, where $\, \alpha = an \,$ and $\, \beta= ma \,$
with $n$ and $m$ fixed (see equations (\ref{beta1}) and (\ref{betaa}).
By proving the conditions of Theorem \ref{thm:bunching}, 
we show that this one-parameter Beta family is monotonically bunching 
in $\alpha$. Finally, Section 5 shows that bunching is invariant under
(simultaneous) strictly monotonic transformations, and gives a real-data
example using the generalized Beta distribution (which is a monotone
transform of the usual Beta) to model income distributions and to provide a new
measure of income concentration.


\section{Definitions}
	\label{chap:def}
\setcounter{figure}{0} \setcounter{equation}{0}

In this section, we present a comprehensive set of foundational definitions that contribute to a precise understanding
of the key concepts and terminologies that will be employed throughout the article.

One of the main tools we use are the well-known inverse probability transform or the
\emph{push-forward} of the measure and the operation it induces on
the corresponding density.
\begin{defn}
	Suppose we have a probability measure $P$ on a space ${\cal{X}}$ and a continuous function 
	$\, y: \, {\cal{X}} \rightarrow {\cal{Y}} \,$.
	The push-forward $\widetilde P$ of $P$ is a measure on the space ${\cal{Y}}$ defined as follows: 
\bsenn
\widetilde P(S):=P(y^{-1}(S)),
\esenn
for a (measurable)
set $S \subseteq {\cal{Y}}$.
\end{defn}

\begin{figure}[!ht]
	\center
	\includegraphics[width=2.5in,height=2.5in]{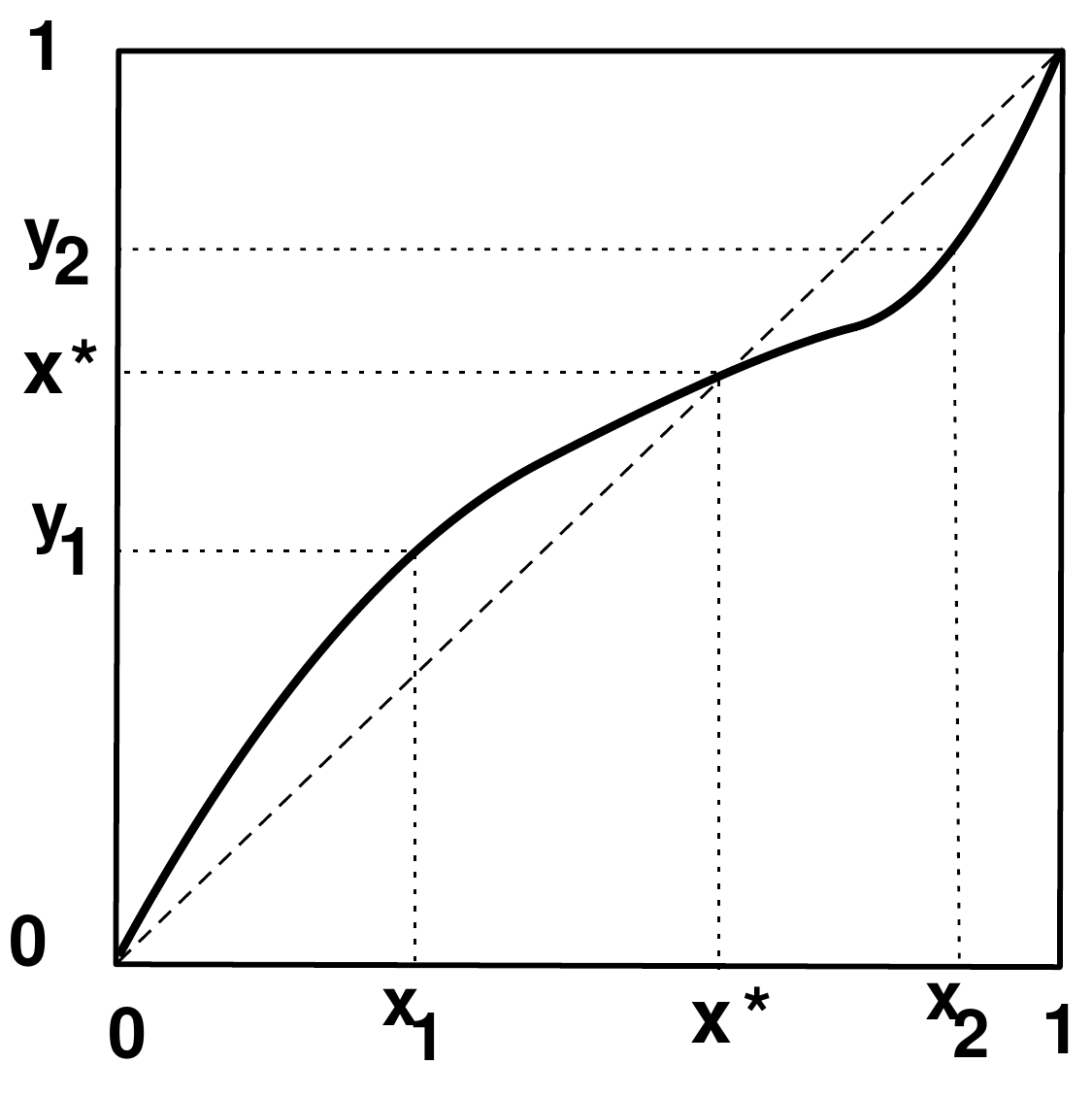}
	\caption{\emph{Definition of the pushforward $\widetilde P$ by $y:X\rightarrow Y$ of a measure $P$.
			The measure of $\widetilde P([y_1,y_2])$ is set equal to that of $P([x_1,x_2])$ where $y_i:=y(x_i)$. }}
	\label{fig:pushforward}
\end{figure}

For the restricted setting: $\, {\cal{X}} = {\cal{Y}} = [0, \, 1] \,$, $\, y: \, {\cal{X}} \rightarrow {\cal{Y}} \,$ invertible,
and with c.d.f. $F$ (on ${\cal{X}}$), the \emph{push-forward} is just the usual probability transform
\bse
\widetilde  F (u) = F(y^{-1}(u)).
\label{eq:def-pushforw}
\ese
Now let's assume in addition that $F$ and $\widetilde F$ and $y$ are continuously differentiable and $y'(x)>0$.
The derivatives of $F$ and $\widetilde F$ are probability densities and will be denoted by $f$ and $\widetilde f$,
respectively, satisfying 
\bse
\widetilde f(y)=\dfrac{f(x)}{y'(x)} \quad \where \quad x=y^{-1}(x) .
\label{eq:perron-frob}
\ese

In this study, we confine our focus to continuous random variables with bounded support within the interval $[0,1]$
such that their probability densities
are contained in the class ${\cal F}_{a}$  defined as following.
\begin{defn}\label{defn:Fa}
	For $a$ in some interval $A\subset\R_{+}$,
	we say $f_a$ (or $F_a$) is in ${\cal F}_{a}$ if 
	$f_a:[0,1]\rightarrow \R$ is a continuous, positive probability density with $f_a(x)>0$ for all $x\in(0,1)$ and $a\in A$. 	
\end{defn}

Next, we present a brief review of some notions of stochastic orders (see \cite{Shak} for an overview of the
different notions of ordering).

\begin{defn}
	Let $X$ and $Y$ be univariate random variables with
	cumulative distribution functions (c.d.f.'s) $F$ and $G$, survival
	functions $\Bar{F}\left(= 1 - F\right)$ and
	$\Bar{G}\left(= 1 - G\right)$.
	\begin{itemize}
		\item[i)] 	$X$ is said to be smaller than $Y$ in the usual stochastic order,
		denoted by $X\leq_{st}Y$, if 
		$E[\phi (X)]\leq E[\phi (Y)]$ for all increasing functions 
		$\phi :	\R\to\R$, provided the expectations exist.
		\item[ii)] 	$X$ is said to be smaller than $Y$ in the convex order,
		denoted by $X\leq_{cx}Y$, if 
			$E[\phi (X)]\leq E[\phi (Y)]$ for all convex functions 
		$\phi :	\R\to\R$, provided the expectations exist.
			\item[iii)] 	$X$ is said to be smaller than $Y$ in the concave order,
		denoted by $X\leq_{cv}Y$, if 
		$E[\phi (X)]\leq E[\phi (Y)]$ for all concave functions 
		$\phi :	\R\to\R$, provided the expectations exist.
		\item[iv)] 	$X$ is said to be smaller than $Y$ in the increasing convex order,
		denoted by $X\leq_{icx}Y$, if 
		$E[\phi (X)]\leq E[\phi (Y)]$ for all increasing convex functions 
		$\phi :	\R\to\R$, provided the expectations exist.
		\item[v)] 	$X$ is said to be smaller than $Y$ in the increasing concave order,
		denoted by $X\leq_{icv}Y$, if 
		$E[\phi (X)]\leq E[\phi (Y)]$ for all increasing concave functions 
		$\phi :	\R\to\R$, provided the expectations exist.
	\end{itemize}
\end{defn}

In broad terms, when $X\leq_{cx}Y$ is satisfied, there is a greater propensity for $Y$ to adopt extreme values
compared to $X$. This implies that $Y$ exhibits higher variability than $X$. Furthermore, it follows 
that $E[X]\leq E[Y]$ and $Var[X]\leq Var[Y]$, given that $Var[Y]<\infty$.
Moreover, $X\leq_{cx}Y$ if, and only if, $X\geq_{cv}Y$, since  if $\phi$ is convex, then $-\phi$ is concave.
On the other hand, if $X\leq _{icx}Y$ then $X$ is both smaller and less
variable than $Y$ in some stochastic sense. It follows that $X\leq _{icx}Y$
implies $E[X]\leq E[Y]$. 
It is clear that $X\leq _{st}Y$ implies $X\leq _{icx}Y$ and $X\leq _{cx}Y$ also implies 
$X\leq _{icx}Y$. In particular, if $E[X]=E[Y]$, then $X\leq _{cx}Y$ if, and
only if, $X\leq _{icx}Y$ (see section 3.A in \cite{Shak}).
Both increasing order relations are related by $X\leq _{icv}Y$
if, and only if, $-Y\leq _{icx} -X$.
It is worth mentioning that the usual stochastic order is known in economics and finance as the first stochastic
dominance (FSD), while the increasing concave order is referred to as the second stochastic dominance (SSD).

Next, we recall a characterization of increasing convex and concave orders based on the number of crossings of
distribution or density functions (see \cite{Shak} or \cite{Lando}).
Let us denote the number of sign changes of a function, $g$, defined on an interval, $I$, with
\[
S^-(g) = S^-(g(x)) = \sup S^-[g(x_1), \ldots, g(x_n)] 
\]
where $S^-[y_1, \ldots, y_n]$ is the number of sign changes of the sequence, $y_1, \ldots, y_n$, where the
zero terms are omitted, and the supremum is extended over all $x_1 < x_2 < \ldots < x_n$ ($x_i \in I$),
$n < \infty$. 

\begin{lem}\label{lem01} 	
	Let $X$ and $Y$ be univariate random variables with
cumulative distribution functions (c.d.f.'s) $F$ and $G$, density
functions $f$ and
$g$, respectively, with finite means.
\begin{itemize}
	\item[i)] If $S^-(F - G) \leq 1$ and the sign sequence starts with $-$, then $X \geq_{icv} Y$ if and only if
	$E(X) \geq E(Y)$, while $Y \geq_{icx} X$ if and only if $E(Y) \geq E(X)$.
	\item[ii)] Let $S^-(f - g) \leq 2$ and the sign sequence begins with $-$. Then, $X \geq_icv{} Y$ if and only if
	$E(X) \geq E(Y)$,
	while $Y \geq_{icx} X$ if and only if
	$E(Y) \geq E(X)$.
\end{itemize}
\end{lem}

\section{Bunching}	\label{chap:general}
\setcounter{figure}{0} \setcounter{equation}{0}

We introduce a novel form of stochastic order based on concentration, which we
call the "Bunching Property". The main aim is to provide reasonable conditions
under which two distributions are ordered according to bunching. We focus on 
distributions in the smooth family ${\cal F}_{a}$ of Definition 2 above, and
provide the following developments. 


\subsection{Basic results}

\begin{prop}\label{prop:pushforward}
	Let $f_{a_1},f_{a_2}\in{\cal F}_{a}$. Then, for all $a_1$ and $a_2$ in $A$, there is a unique
	diffeomorphism $y:(0,1)\rightarrow (0,1)$ such that
	\bsenn
	f_{a_2}(y(x))y'(x)= f_{a_1}(x) \, .
	\esenn
	Furthermore, $y$ can be extended to a continuous function from $[0,1]$ to itself and $y(0)=0$ and $y(1)=1$.	
\end{prop}

\vskip -0.1in\noindent
{\bf Proof} From the definition \eqref{eq:def-pushforw} of the push-forward, we see that $y$ is defined
as the solution of
\bsenn
H(x,y):=\int_0^y\,f_{a_2}(u)\,du-\int_0^x\,f_{a_1}(v)\,dv = F_{a_2}(y)-F_{a_1}(x)=0\,.
\esenn
Since both densities are positive, it follows that $y(0)=0$ and $y(1)=1$.
Now for every given $x$, by continuity (and positivity of both densities), we see that there must be a unique $y$. Furthermore,
since both $\partial_x H$ and $\partial_y H$ are non-zero in $(0,1)^2$, a straightforward application of the
implicit function theorem, tells us that there are differentiable functions $y(x)$ and $x(y)$ so that
$H(x,y(x))=H(x(y),y)=0$. Naturally, these are inverses of one another, and so $y(x)$ is a diffeomorphism
on $(0,1)$. Differentiating $H(x,y(x))$ with respect to $x$ gives the well-known formula in the proposition.
\QED

For the next result, we first need a simple convexity lemma.

\vskip -0.2in\noindent
\begin{lem} Let $g:\R\rightarrow\R$ be twice differentiable on some interval $I$. If
	$x_1<x_2<x_3$ and
	\bsenn
	\frac{g(x_3)-g(x_2)}{x_3-x_2}\leq \frac{g(x_2)-g(x_1)}{x_2-x_1} \,,
	\esenn
	then there is a point $x\in I$ where $g''(x)\leq 0$.
	\label{lem:convex2}
\end{lem}

\vskip -0.1in\noindent
{\bf Proof} The proof is elementary and consists of applying the mean value theorem repeatedly.
\QED


\vskip -0.2in\noindent
\begin{lem}\label{lemma:3fixpts} 
	Let $f_{a_1}$ and $f_{a_2}$ two twice differentiable densities belong to ${\cal F}_{a}$
	with $\, a_2>a_1 \,$ and a push-forward $y$
	as in Proposition \ref{prop:pushforward}. Suppose they satisfy the additional requirements that
	$\lim_{x\searrow 0} y'(x)$ and
	$\lim_{x\nearrow 1}y'(x)$ are greater than one (or tend to infinity) and that
	the second derivative of $f_{a_1}(x)/f_{a_2}(x)$
 is strictly positive on $(0,1)$. Then:
	\begin{itemize}
		\item[i)] $y(x)-x=0$ if and only if $x\in\{0,x^*,1\}$ and
		\item[ii)] $\forall \, x \in(0,x^*)\,:\; y(x)-x>0 \quad \logand \quad \forall \, x \in(x^*,1)\,:\; y(x)-x<0$ and 
		\item[iii)] $f_{a_1}(x^*)\leq f_{a_2}(x^*)$ .
	\end{itemize}	
\end{lem}

\vskip -0.1in\noindent
{\bf Proof} It is sufficient to prove that $y$ has a unique fixed point (i.e. $y(x)=x$) in $(0,1)$. We follow standard
procedure and rewrite the differential equation in Proposition \ref{prop:pushforward} as a 2-dimensional
autonomous system on $(0,1)\times (0,1)$ with reparametrized time:
\bse
\left\{\begin{matrix} \dot x &=& f_{a_2}(y) \\[0.2cm]
	\dot y &=& f_{a_1}(x)   \end{matrix} \right.
\label{eq:ODE3}
\ese
The RHS of this system is continuously differentiable in $(\epsilon,1-\epsilon)\times (\epsilon,1-\epsilon)$ for any
positive $\epsilon$.
We know from
Proposition \ref{prop:pushforward} that there is a unique solution $(x,y(x))$ such that $y(0)=0$ and
$y(1)=1$. By hypothesis, for $x$ close to zero, $y(x)-x>0$, and so $y(x)$ starts out above the diagonal.
Similarly, for $x$ close to 1, $y(x)-x<0$ by the boundary conditions $\lim_{x\nearrow 1}y'(x)$ and
$y(1)=1$. 

We now study the intersections of $\gamma(x):=(x,y(x))$ with the diagonal.
The tangent of the vector field \emph{along the diagonal} $(x,x)$ in the RHS of \eqref{eq:ODE3} is
$t(x):=\frac{f_{a_1}(x)}{f_{a_2}(x)}$. By hypothesis we have
\begin{equation}
t''(x) > 0 \,.
\label{eq:tisconvex}
\end{equation}
Set $f(x):=y(x)-x$. Now let us suppose there are three distinct points $0<x_1<x_2<x_3<1$ in the open
interval $(0,1)$ where $f(x_i)=0$. Then by Lemma
\ref{lem:convex2}, there is a point $v$ where $t''(v)\leq 0$, which contradicts \eqref{eq:tisconvex}.

The conclusion of this reasoning is that there are at most two distinct points where $y(x)=x$.
Since near $x=0$, $\gamma$ is above the diagonal, and near $x=1$ below it, there must be
(by the intermediate value theorem) at least one crossing at $x^*$ from above to below the diagonal.

\begin{figure}[pbth]
	\center
	\includegraphics[height=2.3in]{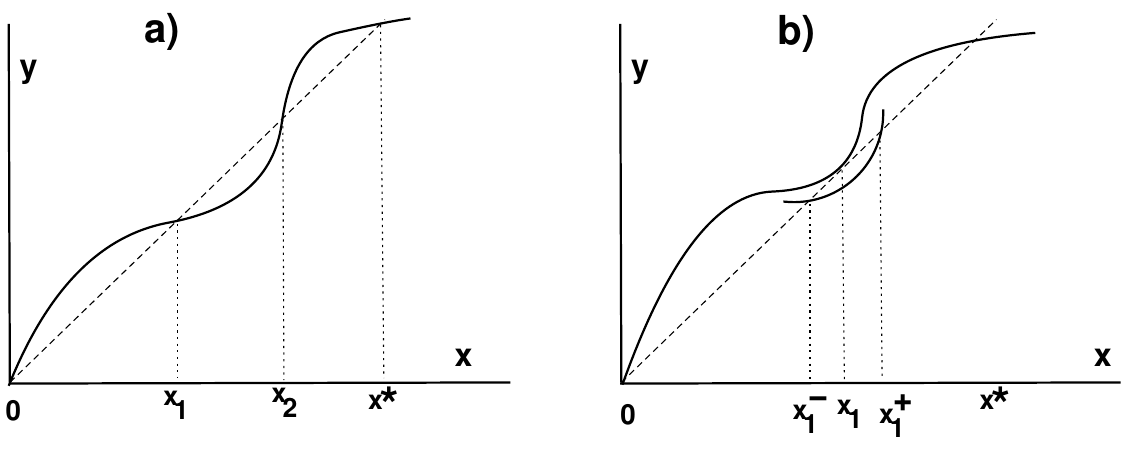}
	\caption{\emph{A non-simple intersection of $y(x)$ and the diagonal at $x_1$. A local twice continuously differentiable
			nearby solution resolves this into two simple intersections.}}
	\label{fig:intersections}
\end{figure}

The only possibilities this leaves is either a unique crossing at $x^*$ or else a crossing at $x^*$ plus
a point $x_1$  where $\gamma$ is tangent to the diagonal. We have to rule out the latter. So suppose
this is the case and assume $x_1<x^*$ as in Figure \ref{fig:intersections}a.
Now solutions to the ODE
of \eqref{eq:ODE3} are unique and cannot cross \cite{arnold,HSD}. Furthermore, nearby solutions approximate one another.
So a solution that starts at $(x_1^-,x_1^+)$ very close to $\gamma$ as illustrated Figure \ref{fig:intersections}b,
must cross the diagonal again at $x_1^+>x_1$ but very close to it. This again
gives three distinct points $0<x_1^-<x_1^+<x^*<1$ such that, respectively,
\bsenn
t(x_1^-)\leq 1\;, \quad t(x_1^+)\geq 1\;, \quad t(x^*)\leq 1 \,.
\esenn
Lemma \ref{lem:convex2} yields a value $x^{**}$ where $t''(x^{**})\leq 0$.
Thus the intersection of $\gamma$ with the diagonal in $(0,1)$ is unique and this gives us items i and
ii of the proposition. Item (iii) follows from the fact that $y(x)$ is differentiable and crosses the
diagonal in the downward direction. \QED

Next, we present a formal statement of the Bunching Property.
\begin{theo}[The Bunching Property]
	Fix $0<a_1<a_2$ and let $F_{a_1},F_{a_2}\in{\cal F}_{a}$ such that they satisfy the conditions of
	Lemma \ref{lemma:3fixpts}.	
 There is a unique $x^*\in(0,1)$ so
	that for all $x_1$ and $x_2$ with $0<x_1<x^*<x_2<1$, we have $ F_{a_1}(x_1)>F_{a_2}(x_1)$ and
	$1-F_{a_1}(x_2) > 1 - F_{a_2}(x_2)$. Thus, we may say that $F_{a_2}$ is more \emph{bunched} around $x^*$
	than is $F_{a_1}$. If $ a_1 > a_2 > 0 $ then $F_{a_2}$ is less \emph{bunched} around $x^*$
         than is $F_{a_1}$.
	\label{thm:bunching}
\end{theo}

\vskip -0.0in\noindent
{\bf Proof} The Theorem is a direct statement of the previous development for $\, a_1 < a_2 \,$. If 
$\, a_1 > a_2 \,$, the same proof goes through with obvious modifications.
\QED

\begin{remark}\label{rem01}
	Observe that, from Theorem \ref{thm:bunching}, we get $S^{-}(F_{a_2}-F_{a_1})=1$ and the sign sequence starts with $-$.
	Therefore, if 
	$X_{a_1}$ and $X_{a_2}$ be two random variables with 
	distribution functions $F_{a_1}$ and $F_{a_2}$, respectively,	
 from Lemma \ref{lem01}(i), we obtain that 
	$X_{a_2}\geq_{icv}X_{a_1}$ if and only if $E(X_{a_2})\geq E(X_{a_1})$,
	while 	$X_{a_1}\geq_{icx}X_{a_2}$ if and only if $E(X_{a_1})\geq E(X_{a_2})$.
\end{remark}

\subsection{The Continuity Property}\label{chap:continuity}

In this subsection, we present a continuity property for the point $x^{*}$ where one distribution is more bunched than the other.

\vskip 0.0in\noindent
\begin{lem} Let $f$ and $g$ be two probability densities that are positive on $(0,1)$. Set
	$\, t(x):=\frac{f(x)}{g(x)}$ and suppose that $t$ is positive and strictly convex with
	\bsenn
	\lim_{x\rightarrow 0}t(x)=\lim_{x\rightarrow 1}t(x)=\infty \,.
	\esenn
	Then $f(x)=g(x)$ in exactly two distinct points in $(0,1)$.
	\label{lem:vu(x)=mu(x)}
\end{lem}

\vskip -0.0in\noindent
{\bf Proof} By the strict convexity, it is clear that $t(x)$ intersects the line $y=1$ either twice or not at all 
But in the latter case, $f(x)$ is strictly larger
than $g(x)$, which conflicts with the fact that both are probability densities (integrating to 1).
\QED


\begin{prop} Let $F_n$ and $G_n$ be c.d.f.'s with densities $f_n$ and $g_n$ continuous in $\, n \,$ and satisfying
	Lemma \ref{lem:vu(x)=mu(x)}. Assume the conditions for Bunching (Theorem \ref{thm:bunching}) and
	let $x^*$ be the unique solution of
	\bsenn
	F_n(x) - G_n(x) = \int_0^x (f_n(t)-g_n(t))\,dt=0 \,.
	\esenn
	Then $x^*$ depends continuously on $n$.
	\label{prop:continuity}
\end{prop}

\vskip -0.0in\noindent
{\bf Proof} 
We first show $\, f_n(x^*) \neq g_n(x^*) \,$. Suppose otherwise. Then there must be a $\, 0<x_-<x^* \,$ such 
that $f_n(x_-)=g_n(x_-)$, for
otherwise one density would dominates the other, and so the integrals could not be equal.
Similarly, since $\, 1 - F_n(x^*) = 1 - G_n(x^*) \,$, the densities must be equal at some point
$\, x_+ \in (x^*, \, 1)$. This would provide three points where the densities are equal, contradicting
Lemma \ref {lem:vu(x)=mu(x)}. Thus $\, f_n(x^*) \neq g_n(x^*) \,$.
Now (by continuity of $f_n$ and $g_n$ in $\, n \,$), a small change in $n$ will cause a small change in the integral:
\bsenn
\int_0^{x^*} (f_{n'}(t)-g_{n'}(t))\,dt=\delta \,.
\esenn
Since near $x^*$ and for $n'$ close enough to $n$, the measures are not equal, a small adjustment
in $x^*$ will bring the integral back to zero.
\QED

\subsection{A Monotonicity Property}
\label{chap:monotonicity}

To provide the main result of this subsection, consider a family of distribution functions 
$\, \{ F_{a;n,m} \} \subset {\cal F} \,$  (see Definition \ref{defn:Fa} ) depending on two additional real parameters: 
$\,  n \geq m > 0 \,$.

\vskip 0.0in\noindent
\begin{cond} Given $0<a_1<a_2$ and $n\geq m>0$. Let $F_a(x)$ and its derivative $f_a(x)$
	be c.d.f.'s and densities as described above. Furthermore, for $\, x \in (0, \, 1) \,$: \\[0.2cm]
	i: If $n=m$, $\, f_a(x) \,$ is symmetric under $x\leftrightarrow 1-x$ and \\[0.2cm]
	ii: $\partial_n f_a(x) \,$ is negative and increasing on $(0,1)$ and \\[0.2cm]
	iii: $\partial_a\partial_n f_a(x)<0 \,$ .
	\label{cond:P-props}
\end{cond}

\vskip -0.0in\noindent
\begin{lem} (See also Proposition \ref{prop:continuity}.) Assume the hypotheses for Bunching 
(Theorem \ref{thm:bunching}) and for Lemma \ref{lem:vu(x)=mu(x)},
Let $x^*$ in $(0,1)$ be the unique point where $y(x^*)=x^*$. Then
	\bsenn
	f_{a_2}(x^*)> f_{a_1}(x^*) \, .
	\esenn
	\label{lem:nu2>nu1}
\end{lem}

\begin{figure}[pbth]
	\center
	\includegraphics[height=4in]{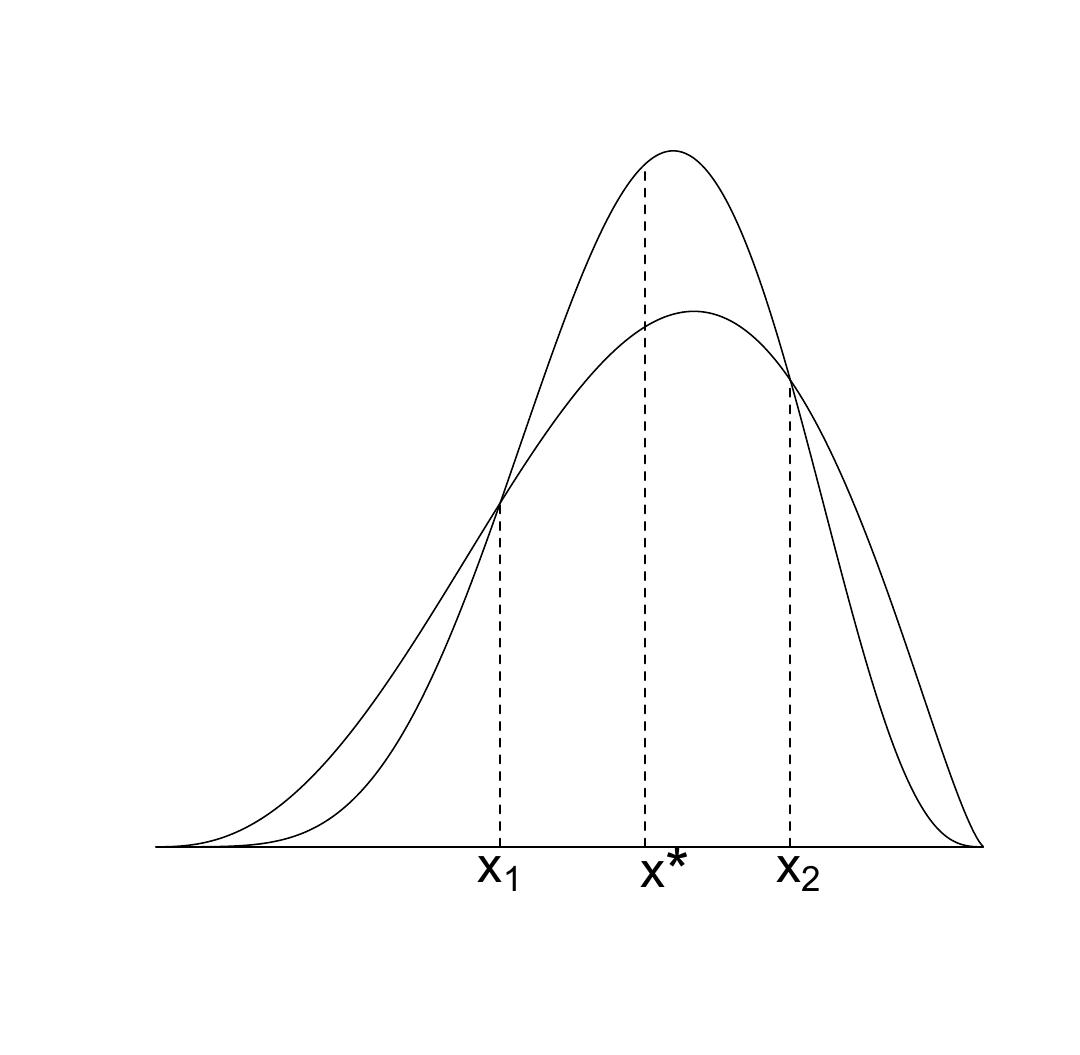}
	\caption{\emph{Two Beta densities ordered by Bunching. Under conditions, 
			if $a_2>a_1>0$, $f_{a_1}$ and $f_{a_2}$
			intersect in exactly two points $\, \{ x_1 , \, x_2 \} \subset (0,1)$.
			At $x^*$, we have that $\int_0^{x^*} f_{a_1}\,dt=\int_0^{x^*}  f_{a_2}\,dt$ .}}
	\label{fig:intersections2}
\end{figure}

\vskip -0.3in\noindent
{\bf Proof.} At the point where $y(x^*)=x^*$, we have $F_{a_2}(x^*)-F_{a_1}(x^*) =0$.
By Lemma \ref{lemma:3fixpts}, $f_{a_1}(x^*)\leq f_{a_2}(x^*)$ and so we only have to rule out
the possibility that they are equal. By Lemma \ref{lem:vu(x)=mu(x)}, $f_{a_2}(x)=f_{a_1}(x)$ in
exactly two distinct points in $(0,1)$, say $x_1$ and $x_2$ as in Figure \ref{fig:intersections2}. 
But since the $f_{a_i}$ are densities (i.e. they integrate
to 1), it impossible for $F_{a_2}(x)$ and $F_{a_1}(x)$ to be equal at either $x_1$ or $x_2$. \QED

In what follows, we will always assume that $a_2> a_1>0$ and $n>m>0$.
From the previous section, we know that there is unique solution for $x$ of
\bsenn
F_{a_2}(x)-F_{a_1}(x)=0 \,.
\esenn
We will track this solution as function of $n$ and $x$, holding fixed all
the other parameters. To facilitate this, we define
\bsenn
J(n,x):=F_{a_2}(x)-F_{a_1}(x) \,.
\esenn

\vskip -0.1in\noindent
\begin{prop} Assume Condition \ref{cond:P-props} and the hypotheses for Bunching 
(Theorem \ref{thm:bunching}) and for Lemma \ref{lem:vu(x)=mu(x)}.
Fix $0<a_1<a_2$ and $n\geq m>0$. We vary $n$ and hold $m$ and the $a_i$ constant.
	For $n=m$, the locus $x^*(n)$ of the unique zero of $J(n,x)$ equals 1/2. We have that
	$x^*$ is a differentiable function satisfying $x^*(m)=1/2$ and $\partial_nx^*>0$.
	\label{prop:derivpositive}
\end{prop}

\vskip -0.0in\noindent
{\bf Proof.} By the symmetry of $f$ (see Condition \ref{cond:P-props}, (ii)) we have that if $n=m$, then $x(n)=1/2$. We also have
\bse
\frac{d}{dn}J(n,x)=\partial_x J\,\frac{dx}{dn} + \partial_n J =0\,.
\label{eq:xprime}
\ese
Now,
\bsenn
\partial_x J(n,x)= f_{a_2}(x)- f_{a_1}(x) \,.
\esenn
By Lemma \ref{lem:nu2>nu1}, have $f_{a_2}(x^*)>f_{a_1}(x^*)$ or
\bsenn
\partial_x J>0 \,.
\esenn
The implicit function theorem now says that near this point there is differentiable function $h$
such that $J(n,h(n))=0$. By Condition \ref{cond:P-props}(iii),we also have
\bse
\partial_n J = \int_0^x \left(\partial_n f_{a_2}(t) - \partial_n f_{a_1}(t) \right) \,dt <0 \,.
\label{eq:TheIntegral}
\ese
Together with \eqref{eq:xprime}, this establishes that the function $x^*(n)$ is differentiable and
strictly increasing, as required. \QED



\section{Applications to Beta distributions}	\label{chap:beta}
\setcounter{figure}{0} \setcounter{equation}{0}

In this section, we present an application of the main results to a beta distribution defined in the interval $[0,1]$. 
To do this, let us recall the definition of this distribution.
 Here we denote the density function of a Beta distribution  with parameters $\alpha>0$ and $\beta>0$ as
\begin{equation} \label{beta1}
f_{\alpha,\beta}(x) = \frac{ x^{\alpha-1}(1-x)^{\beta-1}} {B(\alpha,\beta)} \quad \with \quad B(\alpha,\beta)= 
\int_0^1\,f_{\alpha,\beta}(x) \, \,dx .
\end{equation}
In what follows, we will focus on Beta distributions with parameters
$ n a $ and
$ m a $ for fixed $ n $ and $ m $ with $ a > 0 $. That is, the density and distribution functions of a 
random variable $X\sim Beta(na, ma)$ are denoted:
\begin{equation} \label{betaa}
	p_a(x) = f_{n a,ma}(x) \quad \textup{and} \quad F_a(x) = \int_0^x p_a(t) \, dt \,\, ,
	\end{equation}
respectively. Note that the notation suppresses dependence on $ n $ and $ m $; but when needed we may write
$\, p_{a;n,m}(x) \,$ and $\, F_{a;n,m}(x)$.
From \eqref{betaa}, we see that
\begin{equation*}
	E[X]= \frac{n}{n + m} 
	\quad\textup{and}\quad
	V[X]= \frac{nm}{(n + m)^2(na + ma + 1)}.
\end{equation*}
So the mean is constant and the variance decreases with $a$. This would seem to imply
that the distribution becomes more and more concentrated around the mean.

Next, we prove that the Beta distributions of equation \eqref{betaa} satisfy
the requirements of Proposition \ref{prop:pushforward} and Lemma \ref{lemma:3fixpts}.
First, we need a little lemma.

\vskip 0.0in\noindent
\begin{lem} Let $t(x)=x^{-p}(1-x)^{-q}$. If $p$ and $q$ are
	both positive, then $t$ is strictly convex on $[0,1]$.
	\label{lem:convex}
\end{lem}

\vskip -0.0in\noindent
{\bf Proof.} A tedious, but straightforward, computation gives
\bsenn
t''(x)= \,\frac{qx^2+p(1-x)^2+(qx-p(1-x))^2}{x^{2+p}(1-x)^{2+q}} \,,
\esenn
and the result follows immediately. \QED

\vskip 0.0in\noindent
{\bf Remark.} It is easy to see that the product of two monotone increasing (or decreasing)
convex functions is again convex. In our case, $t$ is the product of one increasing and one decreasing
convex function. It is by no means obvious the the product should be convex. In fact, $p=q=-1$
gives a counter-example.

\vskip 0.0in\noindent
\begin{lem} Let $\{ p_a \} $ be the family of Beta densities (see equation \ref{betaa}) and fix $a_2>a_1$.
	Define $t(x):=\frac{p_{a_1}(x)}{p_{a_2}(x)}$. Then for $y$ in Proposition \ref{prop:pushforward}
	\begin{eqnarray*}
		i: \quad p_a(x)>0 \for x\in(0,1)  \\[.2cm]
		ii: \quad \lim_{x\searrow 0} y'(x) = \lim_{x\nearrow 1}y'(x) = \infty  \\[.2cm]
		iii: \quad  t''(x) >0\; \for x\in(0,1) \\[.2cm]
		iv: \quad \lim_{x\searrow 0} t(x) = \lim_{x\nearrow 1}t(x) = \infty \,.
	\end{eqnarray*}
	\label{lem:inf-deriv}
\end{lem}

\vskip -0.0in\noindent
{\bf Proof.} The first statement is obvious.

The approximate solution for $x$ and $y$ very close to zero can be found by neglecting the $(1-x)$ terms
in the integration (since they are going to be very close to 1). So
\bsenn
G(x,y)\approx \frac{\int_0^y\,u^{na_2-1}\,du}{B(na_2,ma_2)}-
\frac{\int_0^x\,v^{na_1-1}\,dy}{B(na_1,ma_1)} = 0
\quad \Longleftrightarrow \quad \frac{y^{na_2}}{na_2B_2}\approx \frac{x^{na_1}}{na_1B_1}\,.
\esenn
(Here we abbreviated $B(na_i,ma_i)$ as $B_i$.) This gives
\bsenn
y\approx \left(\frac{a_2B_2}{a_1B_1}\right)^{1/na_2}\,x^{a_1/a_2}\,.
\esenn
The first limit follows from $a_2>a_1$. The second limit can be evaluated in the same way
by changing variables $\widetilde x = 1-x$ and $\widetilde y=1-y$. This proves the second statement.

Note that $t(x)=Kx^{-n(a_2-a_1)}(1-x)^{-m(a_2-a_1)}$. The third statement follows from Lemma
\ref{lem:convex}. The fourth statement follows directly from the expression for $t(x)$ just given.
\QED

\begin{theo}
	Let $ \{ F_a(x) \} $ denote the Beta distributions in \eqref{betaa}, 
	and fix $0<a_1<a_2$ and $n\geq m>0$. There is a unique $x^*\in(0,1)$ so
	that for all $x_1$ and $x_2$ with $0<x_1<x^*<x_2<1$, we have $ F_{a_1}(x)>F_{a_2}(x)$ and
	$1-F_{a_1}(x) > 1 - F_{a_2}(x)$. Thus, we may say that $F_{a_2}$ is more \emph{bunched} around $x^*$ than is $F_{a_1}$.
	\label{th_beta}
\end{theo}

\vskip -0.0in\noindent
{\bf Proof.} 
The situation is exactly as sketched in Figure \ref{fig:pushforward} with $p_{a_2}$ on the
vertical axis being the push-forward by $x\rightarrow y(x)$ of $p_{a_1}$. In a case like this,
the bunching property described in Theorem \ref{thm:bunching} holds. This is most easily seen by
noting that the theorem implies that for $y_1\in(0,x^*)$:
\bsenn
F_{a_2}(y_1)= F_{a_1}(x_1)<F_{a_1}(y_1) \quad \where \quad y=y(x) .
\esenn
The latter inequality holds because we also know that $0<x<y(x)<x^*$ and $p_{a_1}(x)>0$ on $[x,y(x)]$.
Similarly, one shows that for $x^*<y_2<1$, $1 - F_{a_2}(y_2) < 1 - F_{a_1}(y_2)$.

\QED

For the family of Beta distributions, the following result is known (see,
e.g., \cite{Arab} or \cite{Lisek}).

\begin{prop} Let $X\sim Beta(\alpha _{1},\beta _{1})$ and $Y\sim Beta(\alpha
	_{2},\beta _{2})$, then $Y\leq _{st}X$ if, and only if, $\alpha _{1}\geq
	\alpha _{2}$ and $\beta _{1}\leq \beta _{2}$.
\end{prop}

In our case, $\alpha _{i}=na_{i}$ and $\beta _{i}=ma_{i}$ for $i=1,2$ and $%
n,m>0$ two real numbers fixed. Therefore, it is evident that $Y\nleq _{st}X$
and $Y\ngeq _{st}X$. 
Note that, from Theorem \ref{th_beta}, we have $S^{-}({F}%
_{a_{2}}-{F}_{a_{1}})=1$ and the sign sequence starts with $-$. Therefore,
from Lemma \ref{lem01}(i), we have the following result.

\begin{corollary} 
Let	$X_{a_{i}}\sim Beta(na_{i},ma_{i})$ for $i=1,2$ and $n,m>0$.
Then, $X_{a_{2}}\geq _{icv}X_{a_{1}}$ for $%
0<a_{1}<a_{2}$. 
\end{corollary}


It is worth mentioning that the previous corollary coincides with Theorem 2 in \cite{lattice}, as $X_{a_i}$
satisfies the conditions of that result for $i=1,2$. Therefore, Theorem 2 in \cite{lattice} can be viewed as a
special case of our Theorem \ref{thm:bunching} for Beta distributions (under certain conditions).

Moreover, since the Beta distributions satisfy all our hypotheses (see Lemma \ref{lem:convex}), the continuity
property holds for them.

\begin{prop}
	Let $ \{ F_a(x) \} $ denote the beta distributions in \eqref{betaa}, 
	and fix $0<a_1<a_2$ and $n\geq m>0$. The location of $x^*$ is a continuous
	function of $n$.
	\label{thm:continuity}
\end{prop}




\section{Application to Income Data with Generalized Beta Fits} \label{sec:data}
\setcounter{figure}{0} \setcounter{equation}{0}

The family of Generalized Beta distributions has been widely studied 
and applied recently, especially for intrinsically positive measurements
like income and lifetime distributions. The 4-parameter Generalized Beta
family is defined as follows: let $s$ be a scale parameter and $\gamma$ a
shape parameter. Then a random variable, $X$, has a Generalized Beta distribution,
$GB2(s, \, \gamma \, \alpha, \beta)$ if
\begin{equation} \label{GB2def}
 \frac{(X/b)^\gamma}{1 + (X/b)^\gamma)} \, \sim \, Beta(\alpha, \, \beta) \, .
\end{equation}
Equivalently,
\begin{equation} \label{GB2def2}
X = b ( Y/(1-Y) )^{1/a} \quad ; \quad Y \sim Beta(\alpha, \, \beta) \, .
\end{equation}

For an introduction to this research, \cite{Chot} provides a general overview, and 
\cite{Spas} is a very recent paper suggesting a new approach to estimation and
references to the recent literature. To apply Theorem \ref{th_beta} to the
GB2 family, the following result is offered:

\begin{theo} Let random variables $X$ and $Y$ have distributions that are ordered
by bunching, and let $T$ be a strictly monotonic function on the real line. Then
$T(X)$ and $T(Y)$ are ordered by bunching in the same direction if $T$ is increasing
and in the opposite direction if $T$ is decreasing.
\label{thm:mono}
\end{theo}

\vskip -0.0in\noindent
{\bf Proof.} 
Bunching is defined by orderings of the c.d.f.'s, and monotone transformations
either maintain ordering (if increasing) or reverse it (if decreasing). 
\QED

Here we consider fitting GB2 distributions to U.S.household incomes
from 2012 to 2022 using Census Bureau data \cite{Guz}.
This data is grouped, giving the percent of households falling
between successive breakpoints 
(0, 15, 25, 35, 50, 75, 100, 150, 200, $>$200) (in units of \$1000). 
While maximum likelihood approaches (among others) have been suggested
for estimating pararmeters from grouped data, we will use minimum 
Chi-Square estimates, which tend to be more numerically stable and
equally efficient in large samples. Minimization was obtained using the
R-program ``optim''.

Of course, the Beta-bunching result, Theorem \ref{th_beta},
applies only if the only varying
parameter is the Beta a-parameter (in (\ref{betaa})). Since the U.S.
household income distribution does not change quickly, it is possible that
the parameters of the fitted distribution are sufficiently close so that
distributions for different years may be compared by bunching. Also, note
that to apply Theorem \ref{th_beta}, estimates of the parameters
$(s, \, \gamma, \, \alpha, \, \beta)$ will not
determine the three Beta parameters $(a, \, n, \, m )$ of (\ref{betaa}). So 
the following parameterization is introduced: $\, a \rightsquigarrow a(n+m) \,$
and $\, \xi \equiv n/(n+m) \,$. Then $\, \alpha = a \, \xi \,$ and
$\, \beta = a(1-\xi) \,$. Thus, estimates for $\alpha$ and $\beta$
uniquely determine estimates for $\xi$ and $\, a = \alpha + \beta \,$, and
bunching order is again monotonic in $a$.

For estimating all 4 GB2 parameters from the income data, the estimates 
were remarkably similar for years 2017 - 2022; but the ``optim''
program seemed to fail in several earlier years: some parameter estimates
reached upper bounds that were imposed, and the minimum Chi-Square
value was much larger than for later years (suggesting that the algorithm
missed the global minimum and was perhaps following a ridge in the Chi-Square
objective function). Thus, to simplify the analysis, the scale parameter
was set to the median income in each year (that is, the breakpoints were
scaled by the median income in each year), and only parameters
$(\gamma, \, \alpha, \, \beta)$ were estimated. These were somewhat constant
over years, and $\gamma$ was generally rather near 1. So to simplify further,
$\gamma$ was set to 1, and the Beta $\xi$ and $\, a \,$ parameters were
estimated.

The $\xi$ estimates ranges from .481 to .484 with a standard deviation of .001;
and so were significantly less than .5 and not substantially non-equal. As a
consequence, if $a_1 > a_2$ denote the (estimated) $a$-values for two years,
then the income distribution for year 1 is less bunched that the distribution
for year 2 (since $\, \xi < .5 \,$ and so $\, \alpha < \beta \,$). Figure
5.1 plots the value of the $a$-estimates over 2012 - 2022, and compares the
a-values to the Gini Index (a classical measure of income inequality).
The Gini-Index values are also taken from Census Bureau etimates (\cite{Kollar}),
where it is noted that the definition of the estimate changes in 2013 and 2017.

Figure 5.1 shows clearly that the Gini Index and bunching can tend in different
directions, with the Gini Index somewhat more erratic. Both show a small trend 
toward more inequality over the period, but the pattern seems clearer with
bunching (which has only two reversals until the final year).  Of course, this
is only a rough informal analysis, and no attempt has been made to provide either
a careful statistical analysis of significance nor a careful economic analysis. 
Nonetheless, it does indicate that bunching should have useful implications for
statistical as well as theoretical questions involving real data applications.

\begin{figure}
	\center
	\includegraphics[height=3in]{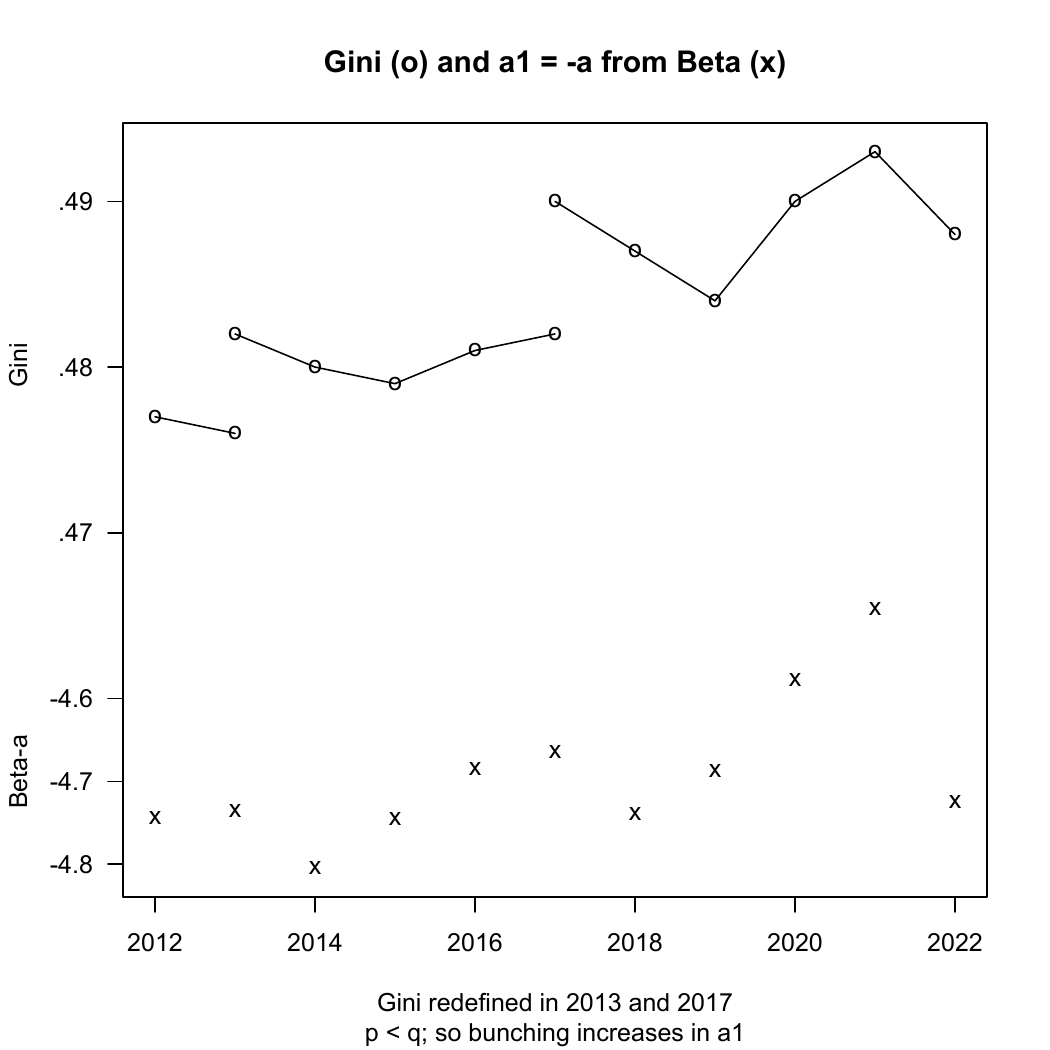}
	\caption{\emph{Comparison of the Gini Index for household
	incomes and the negative of the Beta-a parameter; so increasing
	values of -a also indicate more inequality (concentration)}}
	\label{fig:Gini}
\end{figure}

\section{Conclusions and Motivating Problem}	\label{sec:conclusion}
\setcounter{figure}{0} \setcounter{equation}{0}

We have introduced a novel form of stochastic order based on concentration.
Basically, we say that one distribution is more bunched that another if there is
a point $\, x^* \,$ such that interval probabilities both to the left and right of $\, x^* \,$ 
are smaller (for the more-bunched distribution).
We have provided sufficient conditions for the existence of $\, x^* \,$  and give
conditions for some continuity and monotonicity properties. We also provide an application
to the assessment of income inequality using data
from the Census Bureau on yearly distributions of household income in the U.S. While stochastic
orders are often used in practical statistics to compare distributions generally, this example
suggests that bunching in one-parameter families can create an alternate measure of income
concentration and may offer an useful approach to
comparing certain distributions empirically. 

While these results are of general interest, it may be noted that this research was sparked 
by a specific question from a colleague, Subhash Kochar, who asked if we could
show that $\, P_a({\smhalf} ) \,$ is monotonically decreasing in 
$\, a > 0 \,$ if $\, n > m \,$ for the restricted Beta subfamily.
Apparently, the earliest paper discussing this monotonicity result is a
Clemson University Technical Report, Alam \cite{Alam}, which motivates the problem
in terms of  Ranking and Selection probabilities. The problem also arises
in comparing Gamma distributions in reliability theory: if $\, U \,$ and $\, V \,$
have Gamma distributions with mean parameters $\, n \alpha \,$ and $\, m \alpha \,$
respectively, then $\, \Pr\{U < V \} = P_\alpha(\smhalf) \,$. This suggested
a possible connection to stochastic dominance within the restricted beta family. 

Specifically, if one could prove that $\, x^*(n) > \half \,$ for $\, n > m \,$
in this restricted beta family, then monotonicity of $\, P_a({\smhalf} ) \,$ 
would follow immediately from Bunching for the family (Theorem \ref{th_beta}). Unfortunately, 
it has not been possible to show that Condition \ref{cond:P-props} holds for the Beta subfamily.
It is  possible to show that an alternative sufficient
condition would be $\, \partial_n \partial_a F_{a;n,m}(x) > 0  \,$ for $\, x = x* \,$ and$\,  n > m \,$,
but this also has not been shown for the Beta subfamily. While the desired monotonicity of 
$\, P_a({\smhalf} ) \,$ remains  conjectural, the development of bunching offers a novel and
promising approach to comparing distributions in substantive statistical applications. 

Finally, on generalizing these ideas: clearly our results on Bunching do not 
depend on the specific
domain interval, but do require a one-dimensional (smooth) family in order to 
apply the "push-forward" probability transform. In multivariate situations, 
this can often be replaced by measure transport. See \cite{Hallin} for
a summery of some recent statistical work using this concept. We conjecture that it is 
possible to develop reasonable conditions for smooth multivariate distributions
under which there is a unique point such that the probabilities are monotonically
decreasing for all closed convex sets not containing the point. Such case would 
clearly provide a definition for a strong form of multivariate bunching.

\section*{Funding} 
N. Torrado states that this manuscript is part of the project 
TED2021-129813A-I00 and the
grant CEX2019-000904-S.
She thanks the support of MCIN/AEI/10.13039/501100011033 
and the European Union “NextGenerationEU”/PRTR.


\vspace{\fill}
\end{document}